\documentclass[12pt]{article}
\usepackage{epsfig}
\newcommand{\be}{\begin{equation}}
\newcommand{\ee}{\end{equation}}

\begin{document}
\def\mi{\tilde n}
\def\pe{\frac{\partial}{\partial \Delta E_n}}
\def\nm{n_{\tiny max}}
\def\os{(\frac{\nm-m}{2}+1) }
\def\o{\frac{\nm-m}{2} }
\def\tb{\bar\theta}
\def\l{\langle}
\def\r{\rangle}
\def\nn{\nonumber}
\def\mb{\makebox}
\hspace*{\fill} LMU--TPW--98--11, MPI--PhT--98/62\\
\hspace*{\fill} July 1998\\[3ex]


\def \a{\alpha}
\def \b{\beta}
\def \g{\gamma}
\def \om{\omega}
\def\D{\Delta}
\def\H{\cal H}
\newcommand \reals{I \! \! R}
\newcommand \compl{C \! \! \! \! {\scriptscriptstyle {}^{{}_|}}\ }
\newcommand \N{I \! \! N}
\newcommand \R{I \! \! R}
\newcommand \Z{Z \! \! \! Z}
\newcommand \C{\compl}

\renewcommand{\thefootnote}{\fnsymbol{footnote}}

\begin{center}  
\Large\bf
Convergent Perturbation Theory for a $q$--deformed 
Anharmonic Oscillator \\[4ex]
\normalsize \rm
R.\ Dick$^{a}$, A.\ Pollok-Narayanan$^{a,b}$\footnote{
Andrea.Pollok@physik.uni--muenchen.de}, 
H.\ Steinacker$^{a}$\footnote{Harold.Steinacker@physik.uni--muenchen.de} 
and  J.\ Wess $^{a,b}$ \\[1ex] {\small\it
 ${}^a$Sektion Physik der Ludwig--Maximilians--Universit\"at\\
Theresienstr.\ 37, 80333 M\"unchen,
Germany\\[1ex] 
 ${}^b$Max--Planck--Institut f\"ur Physik \\
 F\"ohringer Ring 6, 80805 M\"unchen, Germany}
\end{center}

\vspace{5ex}
\noindent
{\bf Abstract}: A $q$--deformed anharmonic oscillator is defined 
within the framework of $q$--deformed quantum mechanics.
It is shown that 
the Rayleigh--Schr\"odinger perturbation series for the bounded spectrum
converges to exact eigenstates and eigenvalues, for $q$ close to 1. 
The radius of convergence becomes 
zero in the undeformed limit.
\vspace{\fill}
\noindent

\newpage
\section{Introduction}

The anharmonic oscillator $H=\om a^{\dagger} a + \g X^4$
is a basic quantum mechanical problem with 
one particularly interesting feature: its perturbation series diverges,
but nevertheless there exist eigenstates and energies which are smooth 
as the (positive) coupling constant $\g$ goes to zero \cite{lop,bw,pade}. 
A similar phenomenon is expected to occur in many interacting 
quantum  field theories. The anharmonic oscillator 
can in fact be considered as a $(0+1)$--dimensional 
$\varphi^{4}$ ''field'' theory with  one degree of freedom.

In this paper, we study the analog of this model
in the framework of $q$--deformed quantum mechanics,
based on the $q$--deformed Heisenberg algebra introduced in \cite{pr}.
In particular, one would like to 
know how the perturbation theory of the $q$--deformed anharmonic 
oscillator behaves compared to the undeformed case.
This is of interest in view of a possible
$q$--deformation of  field theory, which 
is expected to be less singular than field
theory based on ordinary manifolds, since $q$--deformation generically
puts physics on a $q$--lattice \cite{pr,bianca}. 
With this motivation, we study the perturbation theory of
the anharmonic oscillator in terms of the $q$--deformed harmonic
oscillator, which was introduced in \cite{macfarlane,biedenharn} and 
realized in the framework of $q$--deformed quantum mechanics in \cite{ho}.

There is considerable freedom in
defining  a $q$--deformed anharmonic oscillator for $q \neq 1$.
Taking advantage of this freedom, we show that for a suitable definition
of the anharmonic oscillator, the perturbation series 
converges to exact eigenvalues and eigenstates for $1<q<1.06$ 
with a certain radius of convergence in $\g$.  
In the limit $q\rightarrow 1$, the model reduces to
the usual anharmonic oscillator, and the radius of convergence goes to zero.
The upper limit on $q$ is not significant.

This paper is organized as follows:
In section 2 we review the $q$-deformed harmonic oscillator 
and its spectrum, and calculate the relevant matrix elements. 
In section 3, the perturbation series for eigenvalues 
and eigenstates is discussed. Some estimates for
the matrix elements are given in the Appendix.

\section{The q-deformed harmonic oscillator}
In this section, we give a brief review of the $q$--deformed
harmonic oscillator, 
and its realization in terms of a $q$--deformed Heisenberg 
algebra. For a more detailed discussion, see \cite{ho} and \cite{pr}.

The q-deformed Heisenberg algebra is the star--algebra 
generated by $X,P,U$ with the relations
\cite{pr}
 \begin{eqnarray}
q^{\frac{1}{2}}XP-q^{-\frac{1}{2}}PX&=&iU\\
UX=q^{-1}XU, \; \;\;  UP&=&qPU. \nn
\end{eqnarray}
We assume $q>1$ to be real. 
The star structure is such that $X$ and $P$ are hermitian, and $U$ is 
unitary:
\begin{eqnarray}
X=X^\dagger,\;\;\; P=P^\dagger,\;\;\; U^{\dagger}=U^{-1}.
\end{eqnarray}

This algebra has the following (momentum--space) representation \cite{pr}:
\begin{eqnarray} 
P|n,\sigma\r&=&\sigma q^n|n,\sigma\r \nn \\
U|n,\sigma\r&=&|n-1,\sigma\r \nn \\
U^{-1}|n,\sigma\r&=&|n+1,\sigma\r \nn \\
X|n,\sigma\r&=&i\sigma\frac{q^{-n}}{q-q^{-1}}(q^{\frac{1}{2}}|n-1,
\sigma\r-q^{-\frac{1}{2}}|n+1,\sigma\r) \nn\\
\l n,\sigma|m,\sigma'\r&=&\delta_{n,m}\delta_{\sigma,\sigma'} 
\label{rep}
\end{eqnarray}
with $n,m\in\N$ and $ \sigma,\sigma' =\pm 1$. 
The completion of these states defines a Hilbert space $\cal{H}$. 

The two values of $\sigma$ describe positive respectively negative momenta. 
(\ref{rep}) is a star--representation, i.e. the star is implemented as 
the adjoint of an operator, and both $X$ and $P$ have  selfadjoint
extensions. That is the reason for introducing $\sigma$,  see \cite{rp}.

This is a starting point for studying $q$--deformed quantum mechanics
\cite{lor,lor2,tm,pr}. In particular, one can define $q$--deformed creation
and anihilation operators as follows:

\begin{eqnarray}
a&=&\alpha U^{-2M}+\beta U^{-M}P\\
a^\dagger&=&\bar{\alpha}U^{2M}+\bar{\beta}PU^M \nn
\label{creat_ops}
\end{eqnarray}
with $M\in\N$, and $\alpha,\beta\in\C$. 
They satisfy the Biedenharn--Macfarlane algebra 
\cite{macfarlane,biedenharn}:
 \begin{eqnarray}
aa^{\dagger}-q^{-2M}a^{\dagger}a=(1-q^{-2M})\alpha\bar{\alpha}=1
\end{eqnarray}
where we fix $\a=\frac i{\sqrt{1-q^{-2M}}}$.
The occupation number operator is defined as
\begin{eqnarray}
\hat n=a^{\dagger}a&=&\alpha\bar\alpha+\beta\bar\beta P^2+
\alpha\bar\beta(U^M+q^MU^{-M})P.
\end{eqnarray}
Now one can write down the following Hamiltonian, which constitutes
the $q$--deformed harmonic oscillator:
\be
H_0=\om a^{\dagger} a
\ee
The spectrum of $H_0$ acting on $\cal{H}$
consists of a bounded spectrum with eigenvalues
$E_n^{(0)} = \om [n]_M = \om \frac{1-q^{-2nM}}{1-q^{-2M}}$ which is 
$2M$--fold degenerate, and an unbounded spectrum with eigenvalues 
$\om (q^{2mM} E_0^{(0)} + \frac{1-q^{2mM}}{1-q^{-2M}})$. 
The $2M$ ground states of the bounded spectrum are
\begin{eqnarray}
|0\r_{\sigma,\mu}^{(M)}=\sum\limits_{n=-\infty}^\infty c_0\left 
(-\sigma\frac{\alpha}{\beta}\right)^nq^{-\frac{1}{2}(Mn^2+Mn+2\mu n)}
|Mn+\mu,\sigma\r , \nn\\
0\leq\mu < M.\;\; 
\end{eqnarray}\\
The existence of an unbounded spectrum beyond 
$E_{\infty}=\frac{\om}{1-q^{-2M}}$ is clear in view of 
(\ref{creat_ops}), since $P$ is an unbounded operator on $\cal{H}$.
For simplicity, we will only consider $M=1$ from now on, and omit the
labels $\mu$ and $M$. 

So far, $\b$ was arbitrary. Requiring that the
$a, a^{\dagger}$ are smooth for $q \rightarrow 1$ and
become the usual (undeformed) creation and anihilation operators
in the limit, one finds \cite{ho} that
\be
\alpha = \frac{i}{\sqrt{1-q^{-2}}},\;\;\beta=\frac{i}{\sqrt{2m\omega}}\nn
\ee
where $m$ is the mass.
For this choice, $H_0$ can be interpreted as a $q$--deformation of the
usual harmonic oscillator, and this will be understood in the following.
The normalized states of the bounded spectrum are
\begin{eqnarray}
|n\r_{\sigma}&=&\frac{1}{\sqrt{[n]}}(a^{\dagger})^n|0\r_{\sigma},
\end{eqnarray}
where $[n] = \frac{1-q^{-2n}}{1-q^{-2}}$.
We define ${\cal H}^{b,\pm}\subset \cal{H}$ to be the closure of
the space spanned by the $|n\r_{\pm 1}$.
As $q\rightarrow 1$, ${\cal H}^{b,+}$ becomes the Hilbert space of the 
usual harmonic oscillator, 
while the unbounded spectrum disappears at infinity, 
and the support of the states with $\sigma=-1$
goes to $-\infty$ in the momentum representation. 
We will thus concentrate on ${\cal H}^{b,+}$.

The eigenstates of $H_0$ can also be written in terms of the $q$--deformed 
Hermite polynomials, which satisfy (see \cite{rh}):
\begin{eqnarray}
\xi H^{(q)}_n(\xi)&=&\frac{\sqrt{q}q^{2n}}{2}(H^{(q)}_{n+1}(\xi)
+2q^{-2}[n]H^{(q)}_{n-1}(\xi))
\label{hermite}
\end{eqnarray}
Defining $\xi=\sqrt{m\omega}X$, one has
\begin{eqnarray}
|n\r_{\sigma}&=&\frac{1}{\sqrt{2}^n[n]!} H^{(q)}_n(\xi)|0\r_{\sigma} \nn.
\end{eqnarray}

Using these Hermite polynomials, it is straightforward to calculate 
the action of $X$
on an eigenstate $|n\r_{\sigma}$, and it follows in particular that
$X \cdot {\cal H}^{b,+} \subset {\cal H}^{b,+}$. This will be important
for the perturbation theory below.

Now we turn to the anharmonic oscillator.
The undeformed anharmonic oscillator is defined by 
$H=\om a^{\dagger} a + \g X^4$ for $\g>0$, thus 
one might naively take the same expression for $q>1$, and study its 
perturbation theory. The relevant matrix elements can be calculated e.g.
using (\ref{hermite}), and we find the following results \cite{andrea}:
\begin{eqnarray}
\l n|X^4|n\r &=&\left(\frac{1}{2m\omega}\right)^2q^{8n+6}
\left([n+1]([n+2]+q^{-4}[n+1]+q^{-8}[n])\right .\nn\\
&&+\left . q^{-8}[n]([n+1]+q^{-4}[n]+q^{-8}[n-1])\right)\nn\\
\l n+4|X^4|n\r&=&\left(\frac{1}{2m\omega}\right)^2q^{8n+14} 
\sqrt{[n+1][n+2][n+3][n+4]}\nn\\
\l n+2|X^4|n\r&=&\left(\frac{1}{2m\omega}\right)^2q^{8n+12}
\sqrt{[n+1][n+2]}\nn\\
&&([n+3]+q^{-4}[n+2]+q^{-8}[n+1]+q^{-12}[n])
\label{matrix_elements}
\end{eqnarray}
They are independent of $\sigma$ which is suppressed. 
All other nonvanishing 
matrix elements can be obtained from those by hermiticity.

Looking at the powers of $q$ in the matrix elements, 
one quickly finds that the perturbation series diverges even faster 
than in the undeformed case. 

However, it is important to realize that there is no reason for
considering the {\em same expression} for $H$ as in the undeformed case;
the only requirement one has to impose is 
that $H$ should reduce to the usual anharmonic oscillator 
as $q\rightarrow 1$. Therefore we might just as well consider the 
Hamiltonian 
\be
H=H_0+\g H'
\label{anharm_osc}
\ee
with 
\begin{eqnarray}
H'&=&\frac{1}{2}  (X^4 \hat Q^5+\hat Q^5X^4),\nn \quad \mb{where}\\
\hat Q &=& (1-a^{\dagger}a (1-q^{-2})).
\end{eqnarray}
$\hat Q$ satisfies
\be
\hat Q|n\r = q^{-2n}|n\r .
\ee
The matrix elements $\l n|H'|m\r$ can be easily obtained from 
(\ref{matrix_elements}), see Figure 2. As is shown in the
Appendix, they have the following upper bound:
\begin{eqnarray}
\l n|H'|m\r < C(q)&:=&  [3]_{4}[2]_{8}q^{-2 \nm+10}[\nm]^2 
 = \frac {4q^{10} [3]_{4}[2]_{8}}{81(1-q^{-2})^2}
\end{eqnarray}
for $1<q<1.06$, where $\nm=\frac {\ln 3}{2\ln q}$. In view of the
results of the next section, 
we define (\ref{anharm_osc}) to be the $q$--deformed anharmonic oscillator.

\section{Perturbation Expansion}

We will use the standard  Rayleigh-Schr\"odinger perturbation formulas
for the eigenstates and eigenvalues of 
\be
H = H_0 + H_1 = H_0 + \g H'
\ee
in terms of the  unperturbed ones,
$H_0 |n\r = E_n^{(0)} |n\r$:

\begin{eqnarray}
\Delta E_n&=&\sum\limits_{k=0}^{\infty}  E_n^{(k)} (\Delta E_n,\g):=
\sum\limits_{k=0}^{\infty} 
\l n|H_1\left(\frac{1}{E_n^{(0)}-H_0}Q_n\left(H_1-\Delta E_n\right)\right)^k
 |n\r\nn\\.
\label{energy_exp}
\end{eqnarray}
where $Q_n = (1-|n\r\l n|)$, and
\begin{eqnarray}
| E_n\r&=& 
   |n\r + \sum\limits_{k=1}^\infty \sum\limits_{n_1,\dots n_{k}\atop n_r
\not= n} \frac{ |n_1\r  \left( \prod\limits_{j=2}^{k}  
 \l n_{j-1} | H_1-\Delta E_n |n_{j}\r  \right)   
 \l n_{k}| H_1 |n\r}{\prod\limits_{j=1}^{k} 
(E_{n}^{(0)}- E_{n_j}^{(0)} ) }.\nn\\
\label{states_exp}
\end{eqnarray}

Strictly speaking, we are of course dealing with a degenerate problem
(since $\sigma=\pm 1$); however as already explained, $X$ and $\hat{Q}$ 
leave ${\cal H}^{b,+}$ invariant, thus the two values
of $\sigma$ do not interfere, and we can restrict ourselves to the
$\sigma=+1$ sector. This will be understood in the following.
We will show that these series in fact converge to exact eigenvalues
and eigenstates of the q-deformed anharmonic oscillator, for a certain
range of $\g$ which depends on $q$.

\subsection{Energy Levels}

If $\g$ and $\D E_n$ are not real, then $H_1$ is
understood to act on the right in the above formulas,
so that the matrix elements can be continued
analytically in $\g$ and  $\D E_n$.
We show first that the sum in (\ref{energy_exp})
is absolutely convergent for $|\D E_n| < \om/5$ and $|\g| < \g(q)$, 
where $\g(q)>0$ provided $q>1$, see (\ref{gamma}). 
Thus the rhs of (\ref{energy_exp}) is an analytic
function of $\Delta E_n$ and $\g$ in that domain,
which can be solved for $\D E_n$ by the implicit function theorem,
defining an analytic function $\D E_n(\g)$.

To see that the sum in (\ref{energy_exp}) is (absolutely) convergent 
for  a certain range of $\D E_n$ and $\g$, we first write   $E_n^{(m)}$ 
more explicitely:
\begin{eqnarray}
 E_n^{(1)}&=&\l n|H_1|n\r \nn\\
 E_n^{(2)}&=& \sum\limits_{n_1\atop n_1 \not= n} 
 \frac{\l n|H_1|n_1\r 
 \l n_1|H_1|n\r }{(E_n^{(0)} - E_{n_1}^{(0)})}  \nn\\
 E_n^{(k)}(\D E_n,\g)&=& 
 \sum\limits_{n_1,n_2,\ldots, n_{k-1}\atop n_r \not= n} 
\frac{\l n|H_1|n_1\r   
 \left(\prod\limits_{j=2}^{k-1} \l n_{j-1}|H_1-\D E_n|n_j\r 
 \right) \l n_{k-1}|H_1|n\r } 
 {(E_n^{(0)} - E_{n_{k-1}}^{(0)}) (E_n^{(0)} - E_{n_{k-2}}^{(0)}) 
\ldots (E_n^{(0)} - E_{n_1}^{(0)}) }  \nn\\
&&\phantom{njghjhkhkjkkkkkkjklkldjfkfjlfjdl}\mb{for}\; k\geq 3  
\label{energy_explicit}
\end{eqnarray}
As is shown in Appendix A, 
the following estimate is valid for $q \in ]1;1.06[$:
\begin{eqnarray}
|E_n^{(k)}(\Delta E_n,\g)|&<&\bar E_n^{(k)}(\Delta E_n,\g,q) := 
\frac{(|\g| C(q))^2
(|\g| C(q) +|\D E_n|)^{k-2} 5^{k-1}}{ ([2]\om q^{-2n})^{k-1}} \nn\\
&&\phantom{ndkfjioejdlfjdjhkjhkjlfdjljkfhgghjdl}\mb{for}\; k\geq 2, 
\label{E_n_estimate}
\end{eqnarray}

The factor 5 comes from the fact that for any given $n_j$, there are only
5 possible $n_{j+1}$ such that the matrix elements in the 
perturbation expansion do not vanish (see (\ref{matrix_elements})).

The series (\ref{energy_exp}) is absolutely convergent if  the following 
condition holds:
\begin{eqnarray}
\left|\frac{ \bar E_n^{(k+1)}}{ \bar E_n^{(k)}}\right|&<&
\theta\;\;\;\mb{for some}\;\;\;\theta<1\nn
\end{eqnarray}
Now 
\begin{eqnarray}
\big|\frac{\bar E_n^{(k+1)}}{\bar E_n^{(k)}}\big|&=&
 5\frac{|\g| C(q) + |\D E_n|}{[2] \omega q^{-2n}}
<\theta\nn,
\end{eqnarray}
and we find that the condition holds e.g. for $|\D E_n| < \om/5$ and
\begin{eqnarray} 
|\g| \leq \g(q) :=\frac{\om([2] q^{-2n}-1)}
 {5 C(q)}.
\label{gamma}
\end{eqnarray}

Therefore we have shown that in this domain, the rhs of (\ref{energy_exp})
defines an analytic function in $\D E_n$ and $\g$.  
Notice that $\g(q)\rightarrow 0$ as $q\rightarrow 1$.

Now consider the equation
\begin{eqnarray}
G(\Delta E_n,\g)&:=&\sum\limits_{k=0}^{\infty} E_n^{(k)}
(\Delta E_n,\g)-\Delta E_n=0. \nn
\label{impl_eq}
\end{eqnarray}  
In the above domain,
this is a uniformly convergent series of analytic functions
(for fixed $q$ in the interval  $]1,1.06[$, say). 
But then using (\ref{energy_explicit}), one sees that
\newline
 $\pe \sum_{k=0}^{\infty} E_n^{(k)}(\Delta E_n,\g)
\Big|_{\D E_{n}=0 \atop \g=0} = 0$, i.e.
\begin{eqnarray}
 \pe  G(\Delta E_n,\g) \neq 0 
\end{eqnarray}
for $\g$ and $\D E_n$ in a neighborhood of 0, by analyticity.
Now the implicit function theorem states that there is a 
function $\D E_n(\g)$ which solves (\ref{impl_eq}) and satisfies 
$\D E_n(0)=0$. Moreover, $\D E_n(\g)$ is analytic in a neighborhood of $0$,
and $|\D E_n|<\om/5$ holds automatically if $\g$ is small enough.


\begin{figure}
\begin{center}
\leavevmode
\epsfysize=1.9in 
\epsfysize=3in 
\epsfig{figure=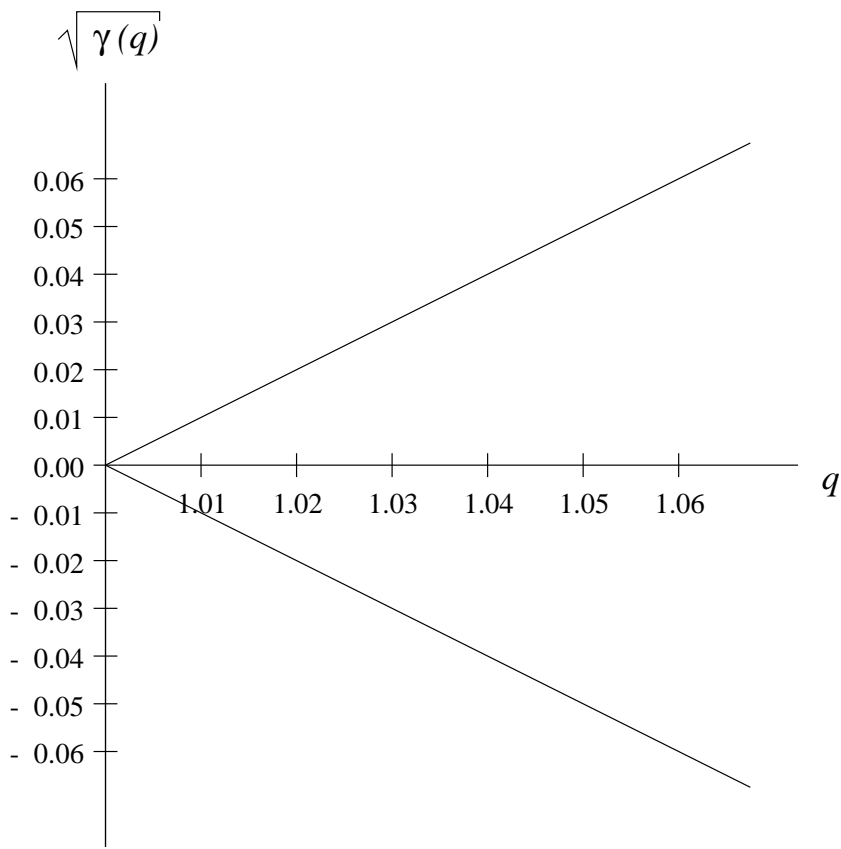,angle=-90}
\caption{Domain of convergence, for $\om=1$}
\end{center}
\label{fig:neu}
\end{figure}

The domain of convergence $\g(q)$  is shown in 
Figure 1 for  $q \in ]1;1.06[$ and $\om=1$. In particular,
$\g(q)$ goes to zero for $q\rightarrow 1 $, in accordance with the 
well--known fact that the perturbation series for the undeformed
anharmonic oscillator is divergent \cite{bw}.

\subsection{Eigenstates}

In this section, we show that (\ref{states_exp}) converges
in ${\H}^{b,+}\subset \H$ for $|\g| <\g(q)$ and $1<q<1.06$, where 
$\D E_n = \D E_n(\g)$ is now the perturbed energy found in the previous 
section. To do this, we have to show that 
\be
\sum_{m=0}^{\infty} |\l m |E_n\r|^2 < \infty,
\ee
or more explicitely

\begin{eqnarray*}
\l E_n|E_n\r &=&\sum\limits_{m=0}^\infty|\l m|E_n\r|^2\nn\\
&=&\sum\limits_{m=0}^\infty \Bigg| \delta_{m,n}+ 
\sum\limits_{k=1}^\infty \sum\limits_{n_1,\dots n_{k}\atop n_r\not= n} 
 \frac{ \delta_{m,n_1} \left( \prod\limits_{j=2}^{k}  
 \l n_{j-1} | H_1-\Delta E_n | n_{j}\r  \right)   
 \l n_{k}| H_1 | n\r}{\prod\limits_{j=1}^{k}(E_{n}^{(0)}- E_{n_j}^{(0)})}
\Bigg|^2\nn
\end{eqnarray*}
From the form of the matrix elements 
(\ref{matrix_elements}), we see that the second term is nonzero only for
$k\geq \frac{|m-n|}{4}$, therefore
\begin{eqnarray}
\l E_n|E_n\r &\leq & 1+ \sum\limits_{m=0}^\infty 
\sum\limits_{k\geq\frac{|m-n|}{4}}  \left( \frac{ (|\g| C(q)+ 
 |\D E_n|)^{k-1} |\g| C(q) 5^{k}}{ ([2]q^{-2n}\omega)^{k} }   
  \right)^2 \nn\\
&\leq& 1+\sum\limits_{m=0}^\infty \left( 5\frac{|\g| C(q) +|\D E_n|}
{[2]q^{-2n}\omega}\right)^{\frac{|m-n|}{2}}
\left(\sum\limits_{k=0}^\infty \left( 
5\frac{ |\g| C(q) + |\D E_n|}{[2]q^{-2n}\om }\right)^{2k} \right) \nn
\end{eqnarray}
for $1<q<1.06$.  
Clearly this converges for $\g$ in the analyticity domain defined above
(such that $|\D E_n| < \om/5$ as before), 
therefore the series (\ref{states_exp}) converges in ${\H}^{b,+}$.
Finally, both $H_0$ and $H'$  leave ${\H}^{b,+}\subset \H$ invariant and 
are bounded operators on ${\H}^{b,+}$ 
($H'$ is bounded because of (\ref{estimate})
and the fact that $H'$ acting on $|n\r$ has no more that 
5 nonvanishing components in terms of that basis). Now it follows that
$|E_n\r$ and $E_n^{(0)}+\Delta E_n$ are indeed eigenstates and eigenvalues
of the full anharmonic oscillator.

As already mentioned, it is known \cite{bw} that the 
undeformed anharmonic oscillator does have nonperturbative eigenstates
and energies for $\g>0$, which are nevertheless smooth as $\g$ goes to zero
from above. Now the formulas (\ref{energy_exp}) ff.
can be analytically continued in $q$ as well, and one would expect that 
the above domain of analyticity for  $\D E_n$ and $\g$  
can be extended to include $q=1$ and 
positive real axis of $\g$. However, at present we are not able to show
this.

\paragraph{Acknowledgements}
A. P.N. acknowledges with thanks the support from MPI.

\section*{Appendix: Matrix elements}
Because $[n]$ is an increasing function in $n$, 
we have the following estimates:

\begin{eqnarray}
\frac{1}{2} \g(1+q^{-40}) q^{14-2n}[n]^2<&\l n+4|H'|n\r&<\frac{1}{2}  
\gamma (1+q^{-40}) q^{14-2n}[n+4]^2\nn\\\nn\\
\frac{1}{2} \gamma (1+q^{-20})q^{12-2n}[4]_{4}[n]^2<&\l 
n+2|H'|n\r&<\frac{1}{2} \gamma (1+q^{-20})q^{12-2n}[4]_{4}[n+3]^2\nn\\\nn\\
 \gamma q^{-2n+6}[3]_{4}[2]_{8}[n]^2<&\l n|H'|n\r& < \gamma 
q^{-2n+6}[3]_{4}[2]_{8}[n+2]^2\nn\\
\end{eqnarray} 
with
\begin{eqnarray}
[n]_{i}:=\frac{1-q^{-ni}}{1-q^{-i}},\;\;\;[n]=\frac{1-q^{-2n}}{1-q^{-2}},\nn
\end{eqnarray}
See Figure 2 for a plot of $\l n|H'|n \r$.

\begin{figure}
\begin{center}
\leavevmode
\epsfysize=1.9in 
\epsfysize=3in 
\epsfbox{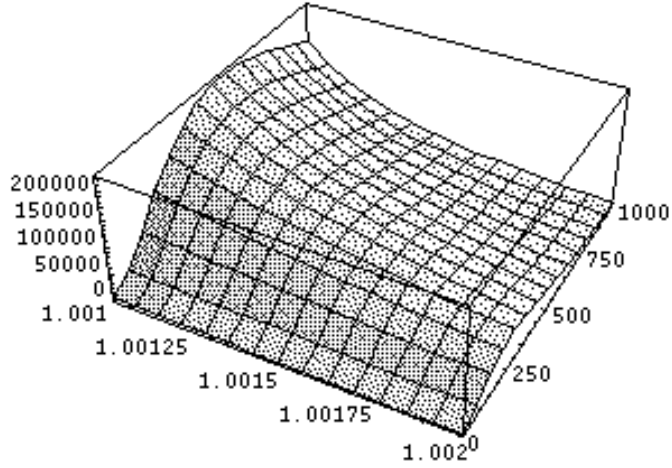}
\caption{The matrix elements $\l n|H'|n\r$ for $q\in [1.001,1.002]$ 
depending on $n$}
\end{center}
\label{fig:matr}
\end{figure}

 To simplify this, 
consider the function $q^{-2n}[n]^2$ for $n\in \reals$, which takes its 
maximum value $\frac {4}{81(1-q^{-2})^2}$ at $n=\nm$, 
\begin{eqnarray}
n_{max}:=\frac{\ln 3}{2\ln q}.
\end{eqnarray}
The matrix elements have a maximum for $n$ close to $\nm$.
More precisely, we can show the following estimate:
\begin{eqnarray}
|\l n+i|H'|n\r| < C(q)&:=& q^{-2\nm+10}[3]_{4}[2]_{8}[\nm]^2
 = \frac {4q^{10} [3]_{4}[2]_{8}}{81(1-q^{-2})^2}
\label{estimate}
\end{eqnarray}
for all $n,m\in\N$. Indeed,
\begin{eqnarray}
\l n+4|H'|n\r&<&\frac{1}{2} (1+q^{-40}) q^{14-2n}[n+4]^2\nn\\
&=&\frac{1}{2} (1+q^{-40})q^{22}q^{-2(n+4)}[n+4]^2\nn\\
&\leq&\frac{1}{2} (1+q^{-40})q^{22}q^{-2\nm}[\nm]^2,\nn
\end{eqnarray}\\
furthermore
\begin{eqnarray}
\l n+2|H'|n\r&<&\frac{1}{2}(1+q^{-20}) q^{18} [4]_4 q^{-2(n+3)}[n+3]^2 \nn\\
    &\leq& \frac 12  (1+q^{-20}) q^{18} [4]_4 q^{-2\nm}[\nm]^2,
\end{eqnarray}
and
\begin{eqnarray}
\l n|H'|n\r &\leq& q^{10}[3]_4 [2]_8 q^{-2\nm} [\nm]^2=C(q)
\end{eqnarray}

Now for  $1\leq q <1.06$, one has
\begin{eqnarray}
1&<&\frac{2 q^{-16}[3]_{4}[2]_{8}}{ 1+q^{-40}} 
\end{eqnarray} (for $i=4$) and
\begin{eqnarray}
1&<&\frac{2[3]_{4}[2]_{8}}{q^{12} [4]_{4}(1+q^{-20})}
\end{eqnarray} 
(for $i=1$).
Combining these estimates, we obtain (\ref{estimate}). 
Furthermore $|E_n^{(0)}-E_{n \pm i}^{(0)}| \geq [i] q^{-2n}\om$, therefore
$|E_{n_j}^{(0)}-E_n^{(0)}| \geq [2] q^{-2n}\om $ in the denominators
of the perturbation expansion, since $i\geq 2$. Now 
(\ref{E_n_estimate}) follows, because for any given $n_j$ 
in the perturbation series, 
there are at most 5 possible $n_{j+1}$ such that 
$\l n_j|H'|n_{j+1}\r$ is nonzero; this means that the number of terms 
at order $k$ is at most $5^{k-1}$.

\end{document}